\documentclass[fleqn,11pt,twoside]{article}

\usepackage{amsthm,amsthm,amssymb, color, xcolor,epsfig, graphics, subfigure}

\usepackage{amsmath, graphicx, latexsym, lscape }

\makeatletter
\newcommand{\copyrightnote}[2]{{\renewcommand{\thefootnote}{}
 \footnotetext{\small\it
\begin{flushleft}
 \copyright \ #1   #2  
\end{flushleft}}}}

\newcommand{\Name}[1]{\begin{flushleft}
                       \LARGE \bf #1
                       \end{flushleft}\vspace{-3mm}}

\newcommand{\Author}[1]{\begin{flushleft}
                       \it #1 \end{flushleft}}

\newcommand{\Address}[1]{\begin{flushleft}
                       \it #1 \end{flushleft}}

\newcommand{\Date}[1]{\begin{flushleft}
                      \small  \it #1 \end{flushleft}}

\newcommand{\evenhead}{Author \ name}
\newcommand{\oddhead}{Article \ name}

\renewcommand{\@evenhead}{
\hspace*{-3pt}\raisebox{-15pt}[\headheight][0pt]{\vbox{\hbox to \textwidth
{\thepage \hfil \evenhead}\vskip4pt \hrule}}}
\renewcommand{\@oddhead}{
\hspace*{-3pt}\raisebox{-15pt}[\headheight][0pt]{\vbox{\hbox to \textwidth
{\oddhead \hfil \thepage}\vskip4pt\hrule}}}
\renewcommand{\@evenfoot}{}
\renewcommand{\@oddfoot}{}

\long\def\@makecaption#1#2{%
  \vskip\abovecaptionskip
  \sbox\@tempboxa{\small \textbf{#1.}\ \ #2}%
  \ifdim \wd\@tempboxa >\hsize
    {\small \textbf{#1.}\ \ #2}\par
  \else
    \global \@minipagefalse
    \hb@xt@\hsize{\hfil\box\@tempboxa\hfil}%
  \fi
  \vskip\belowcaptionskip}

\newcommand{\JNMPnumberwithin}[3][\arabic]{%
  \@ifundefined{c@#2}{\@nocounterr{#2}}{%
    \@ifundefined{c@#3}{\@nocnterr{#3}}{%
      \@addtoreset{#2}{#3}%
      \@xp\xdef\csname the#2\endcsname{%
        \@xp\@nx\csname the#3\endcsname .\@nx#1{#2}}}}%
}

\newcommand{\resetfootnoterule} {
  \renewcommand\footnoterule{%
  \kern-3\p@
  \hrule\@width.4\columnwidth
  \kern2.6\p@}
}

\renewcommand{\footnoterule}{}

\makeatother

\theoremstyle{definition}

\def\a{\alpha}
\def\b{\beta}
\def\G{\Gamma}
\def\A{{\cal A}}

\def\g{{\cal G}}
\def\S{{\cal S}}
\def\M{{\cal M}}
\def\D{\Delta}
\def\be{\begin{equation}}
\def\ee{\end{equation}}
\def\bea{\begin{eqnarray}}
\def\eea{\end{eqnarray}}
\def\ba{\begin{array}}
\def\ea{\end{array}}

\newtheorem{prop}{Proposition}

\begin{document}

\renewcommand{\evenhead}{ {\LARGE\textcolor{blue!10!black!40!green}{{\sf \ \ \ ]ocnmp[}}}\strut\hfill J. Hietarinta and C. Viallet}
\renewcommand{\oddhead}{ {\LARGE\textcolor{blue!10!black!40!green}{{\sf ]ocnmp[}}}\ \ \ \ \  Parametrization of solutions of the Yang--Baxter equations}

\thispagestyle{empty}
\newcommand{\FistPageHead}[3]{
\begin{flushleft}
\raisebox{8mm}[0pt][0pt]
{\footnotesize \sf
\parbox{150mm}{{Open Communications in Nonlinear Mathematical Physics}\ \  \ {\LARGE\textcolor{blue!10!black!40!green}{]ocnmp[}}
\ \ Vol.2 (2022) pp
#2\hfill {\sc #3}}}\vspace{-13mm}
\end{flushleft}}

\FistPageHead{1}{\pageref{firstpage}--\pageref{lastpage}}{ \ \ Article}

\strut\hfill

\strut\hfill

\copyrightnote{The author(s). Distributed under a Creative Commons Attribution 4.0 International License}

\Name{On the parametrization of solutions of\\ the Yang--Baxter equations}

\Author{Jarmo Hietarinta$^{\,1}$ and Claude Viallet$^{\,2}$}

\Address{$^{1}$ Department of Physics and Astronomy, University of Turku\\
        FIN-20500 Turku, FINLAND\\[2mm]
$^{2}$ L.P.T.H.E., Centre National de la Recherche Scientifique \\
         Sorbonne Universit\'e,
        Tour 13, $4^{\rm e}$ \'etage\\
        4 Place Jussieu, F--75252 PARIS Cedex 05,  FRANCE}

\Date{Received Date: October 25,2022; Accepted Date: November 6, 2022}

\setcounter{equation}{0}

\begin{abstract}
\noindent 
We study all five-, six-, and one eight-vertex type
two-state solutions of the Yang-Baxter equations in the form $A_{12}
B_{13} C_{23} = C_{23} B_{13} A_{12}$, and analyze the interplay of
the `gauge' and `inversion' symmetries of these solution. Starting
with algebraic solutions, whose parameters have no specific
interpretation, and then using these symmetries we can construct a
parametrization where we can identify global, color and spectral
parameters.  We show in particular how the distribution of these
parameters may be changed by a change of gauge.
\end{abstract}

\label{firstpage}


\section{Introduction}
The Yang-Baxter equations appeared in the study of two-dimensional
integrable models of statistical mechanics~\cite{Ba81}, and  in the
quantization of 1+1 dimensional integrable equations (see
\cite{Fa82,ZaZa79}). They are an over-determined system of equations
on three matrices $[A,B,C]$ of size $n^2\times n^2$ ($n$ is, e.g, the
number of spin states), and read:
\begin{subequations}\label{eq11}
  \begin{equation}
 \sum_{\alpha_1,\alpha_2,\alpha_3} A^{i_1i_2}_{\a_1\a_2}
                B^{\a_1i_3}_{j_1\a_3}  C^{\a_2\a_3}_{j_2j_3}
     =  \sum_{\b_1,\b_2,\b_3} C^{i_2i_3}_{\b_2\b_3}
          B^{i_1\b_3}_{\b_1j_3}  A^{\b_1\b_2}_{j_1j_2},
\end{equation}
$ \forall i_1,i_2, i_3,j_1,j_2,j_3 = 1\dots n $, or in a shorthand
notation
\be
  A_{12} B_{13} C_{23} = C_{23} B_{13} A_{12}.
\label{E:abcsh}
\ee
\end{subequations}
Here the matrices act on a direct product of three (identical) vector
spaces $V_1\otimes V_2\otimes V_3$, and the subscripts tell on which
spaces the matrix acts non-trivially, e.g., $A_{12}$ means that $A$
acts as $A_{12}\otimes 1$ etc.  We choose to write the equations with
three different matrices $A,B,C$ to emphasize the possible dependence
on parameters, but without prejudice on the nature of these
parameters.

The main issue of the present work is precisely to discuss, by examples,
questions related to the {\em parametrization} of the solutions of
(\ref{E:abcsh}). A natural objective is to write the solution triplet
$[A,B,C]$ in the form of a parametrized family:
\be
A=R(\vec u), \,B=R(\vec u\oplus\vec v),\,C=R(\vec v),
\label{E:parsol}
\ee
using some ``universal'' function $R$.  Here a privileged set of
parameters has been identified, they are the so-called spectral
parameters, and have some kind of addition rule $\oplus$.  The name of
spectral parameter has its origin in the quantum inverse scattering
theory~\cite{Fa82}, and relies on its interpretation as eigenvalue of
a spectral problem.  The spectral parameters play a crucial role in
the quantum inverse scattering approach, and especially in the Bethe
Ansatz construction.  This is the reason why they are singled out in
the parametrization of the solutions.  In Baxter's model \cite{Ba81}
the operation $\oplus$ in (\ref{E:parsol}) is the addition on some
elliptic curve, the uniformization of this curve brings forward
elliptic functions and the moduli of these functions are additional
parameters of the solution: we will call them `moduli parameters'. [It
should be noted that (\ref{E:abcsh}) has many solutions which cannot
be given in the form (\ref{E:parsol}). Take for example $A=B$
arbitrary and $C=P$, the permutation matrix. If this were to be
interpreted according to (\ref{E:parsol}) we should take $\vec v=0$,
$R(0)=P$, and then $R(\vec u)$ remains completely arbitrary.]

In this work we take a closer look on the process by which a good
parametrization can be given to a solution of (\ref{E:abcsh}).  We
show how the inversion symmetries can be used for this purpose.  Of
particular interest is the effect of gauge choice on the nature and
distribution of the parameters.

With reference to the parameter dependence it should be noted that
there are also the so-called `constant Yang--Baxter equations', where
$A=B=C=R$, i.e:
\begin{equation} \label{E:rrr}
  R_{12} R_{13} R_{23} = R_{23} R_{13} R_{12}.
\end{equation}
[For $n=2$ the complete solution of this equation was presented in
\cite{Hi92}.]  Going from (\ref{E:abcsh}) with (\ref{E:parsol}) to
(\ref{E:rrr}), although simple in terms of the parameters (it amounts
to setting them to some value for which $\vec u=\vec v=\vec
u\oplus\vec v$), leads among other things to the successful notion of
the quantum group.  The reverse move, that is to say obtaining solutions
of (\ref{E:abcsh}) starting with solutions of (\ref{E:rrr}), is
sometimes called the `baxterization problem'~\cite{Jo90}, and is
naturally more difficult.  In some cases baxterization is obtained
from group theory~\cite{BeMaVi91d}: There exists a discrete group of
symmetries of equation (\ref{E:abcsh}), the `group of inversions',
which we denote by $\A ut$.  This group acts by non-linear
transformations on the solution triplet and moves it to another
solution. This can be precisely interpreted as the effect of moving
the spectral parameters.  The `baxterization' is essentially the
action of $\A ut$, if it covers densely the manifold of spectral
parameters, but if $\A ut$ produces only a finite set of points we do
not yet have a true baxterization.  The group $\A ut$ is the
statistical mechanical equivalent of the unitarity and crossing
symmetries of $S$-matrix theory: the generators of these symmetries
form a group similar to $\A ut$ \cite{ChLo56,St79,Ba82}.

It is important to note that the symmetry $\A ut$ is not the gauge
symmetry. The latter is believed to bring in only inessential
parameters.  One of the results presented in this paper is that gauge
may also change the distribution of the true parameters (see section
4.2).

The paper is organized as follows: In Sec.\ 2 we discuss in general
the groups of gauge and inversion transformations and their interplay.
In Sec.\ 3 we give all the five- and six-vertex solutions to equation
(\ref{E:abcsh}) and and show how the invariants of the inversion group
can be used to construct a meaningful parametrization for them.  In
Sec.\ 4 we give a seven parameter symmetric eight-vertex solution to
equation (\ref{E:abcsh}) and discuss its parametrization. By allowing
some gauge freedom (for later fixing) we get the solution first in a
rational form.  We then show that a choice of gauge does not change
the nature and number of parameters of the solution triple $[A,B,C]$
but---by its interplay with $\A ut$---affects their distribution
between $A, B$ and $C$, and leads, for example, to Baxter's elliptic
solution.

\section{Two transformation groups}
Suppose $[A,B,C]$ is a triplet of matrices, not necessarily verifying
(\ref{E:abcsh}).  We may define two groups acting on such triplets,
respectively the continuous group of gauge transformations
$\g=SL(n)\otimes SL(n)\otimes SL(n)$, and a discrete group denoted
${\A ut} $.

\subsection{The group of gauge transformations}
Let $g=(g_1,g_2,g_3)$ be an element of $\g$, acting linearly on the
triplet $[A,B,C]$ by similarity transformations:
\begin{equation}
\label{E:gau}
g: [A,B,C] \mapsto
[(g_1\otimes g_2)^{-1} A\, (g_1\otimes g_2),\; (g_1\otimes g_3)^{-1} B\,
(g_1\otimes g_3),\; (g_2\otimes g_3)^{-1} C \,(g_2\otimes g_3)].
\end{equation}
Here the subscript indicates the vector space where the similarity
transformation takes place, and in different spaces the $g$ matrix can
be different.  The group $\g$ is known to take solutions of (\ref{E:abcsh})
into solutions of (\ref{E:abcsh}), but its action is defined everywhere,
even outside the space of solutions.

\subsection{The group of inversions ${\A ut} $}
\label{aut}
Let us first define some elementary operations on a $n^2\times n^2$
matrix $R$, with matrix elements $R^{ij}_{kl}$~\cite{BeMaVi91d}:
\begin{enumerate}
\item{
the (projective) matrix inverse $I$:
\begin{equation}
\sum_{\a\b} (IR)^{ij}_{\a\b} R^{\a\b}_{kl}=\mu \; \delta^{i}_{k} \delta^j_l,
 \qquad i,j,k,l=1,\dots,n,
\end{equation}
with $\mu$ an arbitrary  multiplicative factor.
}
\item{
the transposition $t$:
\begin{equation}
(tR)^{ij}_{kl}= R^{kl}_{ij}, \qquad i,j,k,l=1,\dots,n,
\end{equation}
}
\item{
left and right {\em partial transpositions} $t_l$ and $t_r$:
\begin{equation}
(t_lR)^{ij}_{kl}= R^{kj}_{il}, \qquad (t_rR)^{ij}_{kl}= R^{il}_{kj},
\qquad i,j,k,l=1,\dots,n.
\end{equation}
}
\end{enumerate}
Of course
\begin{equation}
t=t_l\; t_r= t_r\; t_l, \qquad I^2=t^2=t_l^2=t_r^2=1,
\quad \mbox{ and } \quad I\; t = t \; I.
\end{equation}
However,
\begin{equation}
t_l\; I \neq I\; t_l, \quad \mbox{ and } \quad t_r\; I \neq I\; t_r,
\end{equation}
i.e, the two partial transpositions {\em do not} commute with the
inversion, while their product $t$ does.  The transformations $t_l I$
and $t_r I$ are generically of {\em infinite order}, we shall denote by
$\Gamma$ the group generated by $I,t_l,t_r$. [Of course we must assume that
all matrices we are dealing with are nonsingular.]

We may now define the three generators of the `inversion group' as follows:
\begin{eqnarray}
K_a: [A,B,C] &\mapsto & [tI\; A,\qquad t_l\; B,\qquad  t_l\; C], \nonumber\\
K_b: [A,B,C] &\mapsto & [t_l \; A,\quad  t_r I t_l\; B,\qquad
t_r\; C], \label{E:kabc}\\
K_c: [A,B,C] &\mapsto & [t_r\; A,\qquad  t_r\; B,\qquad  tI\; C]. \nonumber
\end{eqnarray}
The three involutions $K_a, K_b, K_c$ act non--linearly (by birational
transformations). They generate an infinite discrete group of
transformations of triplets, which we denote by ${\A ut} $.
\prop
The group $\A ut$ generated by (\ref{E:kabc}) is an invariance group of
the nonsingular solutions of (\ref{E:abcsh}).

\subsection{The compatibility of $\A ut$ with $\g$}
Clearly the action of the two groups $\A ut$ and $\g$ do not commute.
However, their actions are compatible, in the sense that $\A ut$
respects the equivalence classes of triplets $[A,B,C]$ modulo $\g$.
This can be seen as follows: Suppose that $T= [A,B,C]$, and that $T'$
is gauge equivalent to $T$ by $T' = g(T)$, where $g$ acts as defined
in (\ref{E:gau}) with $g=(g_1,g_2,g_3)$. Then from (\ref{E:kabc}) we
get
\[
K_a(T')=  [({^t}g_1 \otimes  {^t}g_2) tI A
({^t}g_1^{-1} \otimes {^t}g_2^{-1}), ({^t}g_1 \otimes   g_3^{-1}) t_l B
({^t}g_1^{-1} \otimes g_3), ({^t}g_2  \otimes  g_3^{-1}) t_l  C
({^t}g_2^{-1}  \otimes  g_3 )]
\]
where ${^t}g_1$ denotes the transpose of $g_1$. This can be written in
the form
\[
K_a \cdot g = g' \cdot K_a \mbox{ with } g'=({^t}g_1^{-1},{^t}g_2^{-1},g_3),
\]
and there are similar relations for $K_b$ and $K_c$. They show that
\prop
If two triplets are gauge related, so are their images by any element of
$\A ut$.
\par\noindent
Note that the previous proposition applies even if the triplet
$[A,B,C]$ does not solve (\ref{E:abcsh}).

\subsection{The moduli space of solutions}
Let $\S$ be any continuously parametrized family of solutions $[A,B,C]$
of (\ref{E:abcsh}).  We will call {\em orbit space} the quotient
\[
\eta=\S/\g
\]
with $\g$ possibly replaced by some of its subgroups. By dividing out
the gauge transformations we obtain the true solution space.
Next we will define the {\em moduli space} of $\S$ by the double
quotient
\[
\M = (\S/\g)/\A ut = \eta/\A ut ,
\]
not caring about the differentiability nor regularity properties of
this quotient. The action of $\A ut$ moves the spectral parameters, so
$\M$ is basically the space of non-spectral parameters.  The second
quotient might be extremely singular. The situation described in this
paper is particularly simple in that respect, since $\eta$ is foliated
by $\A ut$-invariant algebraic subvarieties.

Note that, if $A, B,C$ have a definite form, as is the case for the
five-, six- and eight-vertex Ansatz, we may have to restrict ourselves
to some subgroups of $\g$ and $\A ut$ in order to preserve this form.
We shall in particular need the diagonal subgroup $\g_d$ of $\g$, with
elements
\begin{equation}
g=\left( \left[ \begin{matrix}t_1 & 0 \cr 0 & t_1^{-1} \end{matrix} \right],
         \left[ \begin{matrix}t_2 & 0 \cr 0 & t_2^{-1} \end{matrix} \right],
         \left[ \begin{matrix}t_3 & 0 \cr 0 & t_3^{-1} \end{matrix} \right] \right).
\label{E:diagonal}
\end{equation}

\section{The five- and six-vertex solutions}
\subsection{General considerations}
The six vertex Ansatz for the matrices $A,\,B$ and $C$ is
\begin{equation} X = \begin{pmatrix} 	X_{11} & 0 & 0 & 0 \cr
			0      & X_{22} & X_{23} & 0 \cr
			0  & X_{32} & X_{33} & 0 \cr
			0 & 0 & 0 & X_{44}  \end{pmatrix}	,
\label{E:6veraz}
\end{equation}
For the five-vertex model we take
$X_{32}\equiv 0$, with all other five entries nonzero, while for the
six-vertex model all the six entries are assumed to be nonzero in each
matrix.

The form (\ref{E:6veraz}) is not strictly stable by $\A ut$, since
the partial transpositions exchange the non-zero off-diagonal elements
with the vanishing upper-right and lower-left entries. However, the
subgroup $\A ut_2$ of elements of $\A ut$, which are products of
squares, respects the ansatz.

The action of $\G$ (defined in section (\ref{aut})) on a generic
$4\times 4$ matrix was analyzed in~\cite{BeMaVi92}, where it was shown
that the invariants of $\A ut_2$ are ratios of some quadratic
polynomials in the entries of the matrix.  Out of the 18 polynomials
$p_i$ of~\cite{BeMaVi92}, only five are non-vanishing when evaluated
on a matrix of the form (\ref{E:6veraz}), they are
\begin{eqnarray}
p_1(X) &=& X_{11}X_{22}+X_{33}X_{44}, \nonumber \\
p_2(X) &=& X_{11}X_{22}-X_{33}X_{44}, \nonumber \\
p_5(X) &=& X_{11}X_{33}+X_{22}X_{44}, \\
p_6(X) &=& X_{11}X_{33}-X_{22}X_{44}, \nonumber \\
p_9(X) &=& X_{11}X_{44}+X_{22}X_{33}-X_{23}X_{32}.\nonumber
\end{eqnarray}
Invariants of $\A ut_2$ can then be obtained by taking ratios of the
form $p_i(A)/p_j(A)$, resp.\ $(B), (C)$.  One should notice that,
in the case under study,  all these polynomials are
{\it independent of the gauge parameters},
contrary to what happens for the general (16-vertex) case.

The rank of the system of the four invariant ratios constructed
from the generic matrix (\ref{E:6veraz}) is only 3, the additional
relation being
\begin{eqnarray}
p_1^2-p_2^2 = p_5^2-p_6^2.
\end{eqnarray}
A solution for which $p_9(X)\equiv 0$ is called
``free-fermion type''~\cite{FaWu70,JaMa83b,No92}.

\subsection{The five-vertex solutions}
For the five-vertex solution we take the matrix elements
 $A_{32}=B_{32}=C_{32}=0$.
In addition let us scale so that $A_{11}=B_{11}=C_{11}=1$. Substituting
this Ansatz into (\ref{E:abcsh}) leads easily to precisely two solutions.

\subsubsection{The first solution}
The first solution is given by
\be
 A = \begin{pmatrix}
1&0&0&0\cr
0&x_2&a&0\cr
0&0&x_3&0\cr
0&0&0&x_4 \end{pmatrix},\quad
B = \begin{pmatrix}
1&0&0&0\cr
0&x_2&b&0\cr
0&0&y_3&0\cr
0&0&0&\frac{x_4y_3}{x_3} \end{pmatrix},\quad
C = \begin{pmatrix}
1&0&0&0\cr
0&\frac{x_2}{x_4}&c&0\cr
0&0&y_3&0\cr
0&0&0&\frac{y_3}{x_3} \end{pmatrix},
\label{E:5vs1r}
\ee
with the constraint
\be
ac=b\,(1-x_2x_3/x_4).
\ee

Since now $p_9\ne0$ let us consider
\be
\D:=\frac{p_1^2-p_2^2}{p_9^2},
\label{E:D}
\ee
we find
\be
\D(A)=\D(B)=\D(C)=\frac{4x_4/(x_2x_3)}{[1+x_4/(x_2x_3)]^2}\,,
\ee
and thus we have found a global invariant, which we may exchange for
$d:=\sqrt{x_2x_3/x_4}$ [and in terms of $d$, the constraint becomes
$ac=b(1-d^2)$].  Note that $d$ may be constructed from the covariants
${\cal K}_1,\;{\cal K}_2,\; {\cal K}_3$ of~\cite{BeMaVi92}, while the
modular invariant $J$ of~\cite{BeMaVi92} vanishes.

{}From the other ratios let us look at the following:
\be
\delta := \frac{ p_1+p_2}{ p_1-p_2 },\quad
\delta' := { p_5-p_6}{ p_5 + p_6}.
\label{E:dd}
\ee
We find
\be\begin{array}{rcl}
q_1^2&:=&\delta'(A)=\delta'(B)=(x_2/d)^2,\\
q_2^2&:=&\delta(A)=\delta'(C)=(d/x_3)^2,\\
q_3^2&:=&\delta(B)=\delta(C)=(d/y_3)^2.
\end{array}
\label{E:5v1a}
\ee
Since the common index between $A$ and $B$ is 1, between $A$ and $C$
is 2 and between $B$ and $C$ is 3, we have also introduced new `color'
variables $q_i$ in (\ref{E:5v1a}).  After this let us define
\be
R_{5a}(i,j)=\begin{pmatrix}
1&0&0&0\cr
0&d\,q_i&(1-d^2)\,g_ig_j^{-1}&0\cr
0&0&d\,q_j^{-1}&0\cr
0&0&0&q_i\,q_j^{-1} \end{pmatrix},
\ee
and then the solution (\ref{E:5vs1r}) can be written as
$A=R_{5a}(1,2),\,B=R_{5a}(1,3),\,C=R_{5a}(2,3)$. In this form $g_i$
are the variables that are changed by gauge, we could fix them by
putting $g_i=1$.

Thus, starting with the solution (\ref{E:5vs1r}) without any
particular structure we were able to put it into a form in which the
variables were either global ($d$) or associated with the vector
spaces ($q_i$). [These latter ones are often called `color' variables
\cite{Color}.]  This was accomplished by looking at what remains from
the generic invariants of $\A ut_2$ for the specific solution.  The
 action of $\A ut_2$ on $R_{5a}$ is given as follows:
\be
(K_aK_b)^2: g_1 \mapsto g_1 d^2, \quad
(K_bK_c)^2: g_3 \mapsto g_3 d^2, \quad
(K_cK_a)^2: g_2 \mapsto g_2 d^2
\ee
As a consequence the action of $\A ut_2$  cannot be distinguished from the
one of the gauge transformations, i.e. $\A ut_2$ acts as unity on the
orbit space $\eta$.

At this point let us introduce a notation for the {\em parameter
content} of a solution triple: we say that the parameter content of
$(A,B,C)$ is $(n_A,n_B,n_C)$, if fixing $A$ fixes $n_A$ parameters,
then fixing $B$ fixes $n_B$ of the {\it remaining} parameters and so
on. The parameter content of this solution is clearly $(3,1,0)$, since
from $A$ we get $q_1,q_2,d$ and from $B$ the remaining $q_3$.

Finally we note that in the constant limit we must take all
$q_i$ equal and obtain the well know solution of (\ref{E:rrr})
\be
R=\begin{pmatrix}
1&0&0&0\cr
0&p&1-pq&0\cr
0&0&q&0\cr
0&0&0&1 \end{pmatrix},
\ee
where $p=dq_i,\,q=d/q_i$.

\subsubsection{The second solution (free fermion type)}
By inspection we can write the second solution in terms of
\be
R_{5b}(i,j)=\begin{pmatrix}
1&0&0&0\cr
0&p_i&g_{ij}&0\cr
0&0&q_j&0\cr
0&0&0&-p_iq_j \end{pmatrix},
\label{E:5v2}
\ee
as $A=R_{5b}(1,2),\,B=R_{5b}(1,3),\,C=R_{5b}(2,3)$, together with the
constraint
\be
g_{12}g_{23}=g_{13}(1-p_2q_2).
\label{E:5bg}
\ee
The resolution of this last constraint is more problematic and the
choice of gauge less trivial than in the first solution. [Recall
that under a gauge transformation $g$ (see (\ref{E:diagonal})),
$g_{ij}\to g_{ij}t_i/t_j$, thus we can fix two of $g_{ij}$'s.]

For a uniform representation in which each $R_{5b}(i,j)$ depends only
on two variables we should have $g_{ij}=g(p_i,q_j)$, but it is easy to
see that (\ref{E:5bg}) does not have solutions of this type. We  must
then relax the condition, and if we instead allow $g_{ij}=
g(p_i,p_j,q_i,q_j)$ then a family of solutions can be constructed:
\be
g_{ij}=(1-p_iq_i)^\alpha(1-p_jq_j)^{1-\alpha}.
\label{E:unig}
\ee

The total number of parameters of (\ref{E:5v2}) is four, as in the
previous case, but now they are all `color' parameters
$(p_1,p_2,q_2,q_3)$. If now we choose, for example, the non-uniform
gauge $g_{12}=g_{13}=1$ the parameter content of the solution is
$(2,1,1)$. In the uniform gauge (\ref{E:unig}) with $\alpha=0$ or $1$
we must introduce extra parameters and get parameter content $(3,2,0)$
and $(3,1,1)$, respectively, and in other cases $(4,2,0)$. Thus the
price we have to pay for uniformity is the introduction of extra
spurious parameters, that could in principle be gauged away.

In this case the action of $\A ut_2$ is:
\begin{eqnarray*}
& (K_aK_b)^2 &: g_{12}\mapsto -g_{12},\quad g_{13}\mapsto -g_{13},
\quad g_{23}\mapsto g_{23}, \\
& (K_bK_c)^2 &: g_{12}\mapsto g_{12},\quad g_{13}\mapsto -g_{13},
\quad g_{23}\mapsto -g_{23}, \\
& (K_cK_a)^2 &: g_{12}\mapsto -g_{12},\quad g_{13}\mapsto g_{13},
\quad g_{23}\mapsto -g_{23},
\label{E:fft}
\end{eqnarray*}
which again is indistinguishable from simple gauge transformations of
square one.

The constant limit of this solution is obtained from (\ref{E:5v2}) with
$p_i=p,\,q_i=q,\,g_{ij}=1-pq,\,\forall i,j$.

\subsection{The six-vertex solutions}
In this case all entries in (\ref{E:6veraz}) are nonzero. We use a
gauge transformation with a diagonal $g_i$ and overall scaling to make
$X_{23}=X_{32}=1$ for $A$ and $B$, say. From equation (\ref{E:abcsh})
we then get $C_{23}=C_{32}$, which can be scaled to 1. We will
therefore only write down the diagonal elements as $dp(X):=[X_{11},
X_{22}, X_{33}, X_{44}]$, and simplify notation by using
$x_i:=X_{ii}$, e.g., $dp(A)=[a_1,a_2,a_3,a_4]$.

When the Ansatz (\ref{E:6veraz}) and $X_{23}=X_{32}=1$ is used in
(\ref{E:abcsh}) we get 6 equations,
\[
\begin{array}{rcl}
b_2a_1 - b_1a_2 - c_2&=&0,\\
  - b_1 - c_2a_3 + c_1a_1&=&0,\\
  - a_3 - c_3b_1 + c_1b_3&=&0,\\
 a_2 - c_4b_2 + c_2b_4&=&0,\\
 b_4 - c_4a_4 + c_3a_2&=&0,\\
 b_4a_3 - b_3a_4 + c_3&=&0.
\end{array}
\]
We solve $a_1$ from the second equation, $a_2$ from the fourth, $a_3$
from the third and $a_4$ from the fifth, and then $b_4$ from the
first. Since the parameters are assumed to be nonzero there is no
ambiguity in doing this and what then remains is one equation which
factors as (all computations were done using REDUCE \cite{REDUCE} and
Maple \cite{Maple}):
\be
(b_1 b_2 b_3 + c_3 c_4 b_1^2 b_2 - c_2 c_3 b_1
b_2 b_3 - c_1 c_4 b_1 b_2 b_3 - c_1 c_2 b_3 + c_1 c_2 b_2 b_3^2)
(1 - c_2 c_3 - c_1 c_4)=0.
\label{E:6vfactor}
\ee
We thus recover the two known six-vertex-type solutions:

\subsubsection{The first solution: the asymmetric six-vertex}
\label{asymmetric}

The first solution is obtained when we use the first factor of
(\ref{E:6vfactor}). After solving for $b_2$ and some simple parameter
changes we can write the solution as follows:
\be
\begin{array}{rcl}
dp(C)&=&[a,b,c,d],\\
dp(B)&=&[a e,b f/h,c f,d e/h],\\
dp(A)&=&[ e + b c (f - e),b d (f - e))/h,c a (f - e),(e + b c (f - e))/h],
\end{array}
\label{E:6vsol2}
\ee
where
\be
h:=e f + (a d e-b c f)(e - f).
\ee
In this form the solution has a rational expression, but one cannot
yet identify any spectral (or color) parameters.

Before analyzing the solution in detail let us just note that the
choice of $C$ fixes four parameters, and the further choice of $B$
fixes the two remaining ones, so the parameter content is
(4,2,0).

The group of gauge transformations $\g$ consists of the diagonal
transformations $\g_d$ as above, and they only change the off-diagonal
elements of the matrices. The orbit space $\eta$ is thus six
dimensional and we may take $a,b,c,d,e,f$ as the coordinates on
$\eta$.

In order to construct a good parametrization we start with (\ref{E:D})
and in this case get
\begin{eqnarray}\label{global}
\D(A) = \D(B) = \D(C) = \frac{4 a b c d}{ (a d + b c -1 )^2 },
\end{eqnarray}
which is a `global invariant'.  Here again, the modular invariant $J$
of~\cite{BeMaVi92} vanishes.

{}From (\ref{E:dd}) we find
\begin{eqnarray}
q_3^4  &:=& \delta(C) = \delta(B) = \frac{ab}{cd},  \nonumber \\
q_2^4  &:=&  \delta'(C) = \delta(A) = \frac{bd}{ac}, \\
q_1^4 &:=&  \delta'(B) = \delta'(A) = \frac{bd}{ach^2},\nonumber
\end{eqnarray}
which defines three new (color) parameters. Thus we have been able to
identify four of the six parameters. To study the remaining ones we
note that any matrix of the form (\ref{E:6veraz}) with
$X_{23}=X_{32}=1$ (solution or not) can be parametrized equally well
with the four parameters $\D,\,\delta,\delta',\,\lambda$, (where
$\D,\,\delta,\delta'$ were defined in (\ref{E:D}, \ref{E:dd}) and the
nature of $\lambda$ is left open at the moment) and may be written as
$R(\D, \delta', \delta, \lambda)$. For the present solution we have
\begin{eqnarray}
 A  =  R (\D, q_1, q_2, \lambda_A), \qquad
 B  =  R (\D, q_1, q_3, \lambda_B), \qquad
 C  =  R (\D, q_2, q_3, \lambda_C)
\end{eqnarray}
where $\lambda_A, \lambda_B, \lambda_C$ must verify an additional
relation, which will appear as we clarify the $\lambda$
dependence. For this purpose let us write the matrix elements of the
solution matrix $R(\D, \delta', \delta, \lambda)$ as
\begin{equation}
R_{11}=u\left(\frac{\delta}{\delta'}\right)^{1/4},\quad
R_{22}=v\left({\delta'}{\delta}\right)^{1/4},\quad
R_{33}=v\left(\frac1{\delta'\delta}\right)^{1/4},\quad
R_{44}=u\left(\frac{\delta'}{\delta}\right)^{1/4}.
\end{equation}
This form is compatible with the previous assignments, and the $\lambda$
dependence is entirely inside the $u$'s and the $v$'s. What
remains are the following relations:
\bea
\sqrt\D&=&\frac{2u_Av_A}{u_A^2+v_A^2-1}=\frac{2u_Bv_B}{u_B^2+v_B^2-1}
=\frac{2u_Cv_C}{u_C^2+v_C^2-1},\label{E:cur6v}\\
v_A&=&v_B u_C-v_C u_B,\phantom{|^|}
\quad u_B=u_A u_C-v_A v_C,
\eea
which are resolved by
\bea
u_I&=&\frac{\sin(\gamma-\lambda_I)}{\sin(\gamma)},\quad
v_I=\frac{\sin(\lambda_I)}{\sin(\gamma)},\quad I=A,B,C \\
\lambda_B&=&\lambda_A+\lambda_C,\quad \sqrt\D=-1/\cos(\gamma).
\eea
After changing from $\Delta$ to $\gamma$ we get
\be
\label{E:6V2}
R(\gamma,q',q,\lambda):=\begin{pmatrix}
\dfrac{q}{q'}\,\dfrac{\sin(\gamma-\lambda)}{\sin(\gamma)} & 0 & 0 & 0 \cr
0      &  q q'  \,\dfrac{\sin(\lambda)}{\sin(\gamma)}& 1 & 0 \cr
0  & 1 & \dfrac1{q q'}  \,\dfrac{\sin(\lambda)}{\sin(\gamma)} & 0 \cr
0 & 0 & 0 & \dfrac{q'}{q}\,\dfrac{\sin(\gamma-\lambda)}{\sin(\gamma)} \end{pmatrix}
\ee
and the solution is given by
\be
A=R(\gamma,q_1,q_2,\lambda_A),\quad
B=R(\gamma,q_1,q_3,\lambda_A+\lambda_C),\quad
C=R(\gamma,q_2,q_3,\lambda_C).
\label{E:asym6}
\ee

Thus we have been able to introduce a good parametrization to the
algebraic solution (\ref{E:6vsol2}).  This is the asymmetric
six-vertex solution of~\cite{FaWu70,JaMa83b,No92}.

With this parametrization we see that the spectral parameters of $A$,
$B$, and $C$ are points on the circle (\ref{E:cur6v}), with its simple
addition law.  It also clarifies the action of $\A ut_2$.  Since
\be
 t_l I t_l I : \; R(\gamma, q, q', \lambda)
\longmapsto R(\gamma, q, q', \lambda + 2 \gamma),
\label{E:shift}
\ee
the action of $\A ut_2$ is just a shift of the spectral parameter.
Moreover, among the moduli parameters, only $\gamma$ is global, while
the $q_i$'s are attached to the vector spaces on which the matrix
operates (`color parameters'). It is interesting to note how the
periodic orbits of $\A ut$ appear: they correspond to the values of
$\gamma$ which are commensurate to $\pi$.  The special case
$\gamma=\pi/2$ is included in the following.

The first known parametrized solution of (\ref{E:parsol}),
$R(u)=P+uI$, is obtained as a singular limit of (\ref{E:6V2}): take
$q_k=\sqrt{-1}$, $\lambda=-\gamma u$ and then let $\gamma\to0$. The
case 6V(I) of \cite{SoUcAkWa82} is sub-case $q=q'$ of
(\ref{E:6V2}). For other special cases, see \cite{SuWaWu94}.

\subsubsection{The second solution (free fermion type)}
If we solve for $c_4$ from the second factor of (\ref{E:6vfactor}), we
get $C$ and $B$ as follows:
\be
\begin{array}{rcl}
dp(C)&=&[c_1,c_2,c_3,(1-c_2 c_3)/c_1],\\
dp(B)&=&[b_1,b_2,b_3,(1-b_2 b_3)/b_1],\\
\end{array}
\label{E:61CB}
\ee
In this case it is natural to put the diagonal elements of
(\ref{E:6veraz}) into a $2\times 2$ matrix as
\be
\hat X=\begin{pmatrix}X_{11} & -X_{22}\cr X_{33} & X_{44} \end{pmatrix},
\ee
and then (\ref{E:61CB}) implies $\det \hat B=\det \hat
C=1$. Furthermore, we find that the remaining matrix $\hat A$ is given
by a matrix product (note the order)
\be
\hat A=\hat C^{-1} \hat B.\label{E:slc}
\ee
Thus in this case the natural parametrization is through the group
$SL(2)$: For any $SL(2)$ matrix $\hat X$ let $R(\hat X)$ be the
$4\times4$ matrix of type (\ref{E:6veraz}) obtained by putting the
elements of $\hat X$ on the diagonal as discussed above
($X_{23}=X_{32}=1$).  Then, according to the above, we can write the
result in the form $A=R(\hat A),\,B=R(\hat A\oplus \hat C),\,C=R(\hat
C)$, where now $\hat A\oplus \hat C=\hat C\hat A$ (in particular, here
$\oplus$ is not Abelian). The parameter content of this solution is
clearly $(3,3,0)$. This already shows the difference with the first
solution.  This solution is the one of~\cite{IzKo82}, see
also~\cite{BaSt85a}. Its constant limit is the permutation matrix.

There are no gauge parameters in our presentation of the  solution.
The action of $\A ut_2$ is
\be
 (K_aK_b)^2: 	\hat{A} \mapsto -\hat{A},\quad
 \hat{B}\mapsto -\hat{B},\quad
 \hat{C} \mapsto \hat{C}
\ee
and so on. It is equivalent to discrete gauge transformations of
square one.

This solution allows many reductions with one-dimensional spectral
parameters.  For example, we may take the solution (\ref{E:6V2}), with
$\gamma=\pi/2$.  Consider the polynomials $p_i$ introduced above.  In
this case $p_9$ vanishes, and the rank of the remaining invariant
ratios is 2. Fixing the value of $\delta= q^4$ and $\delta'=q'^4$
determines a curve on $SL(2)$, leading to an Euler type of
parametrization:
\be
\hat X = 
\left [\begin {array}{cc} q&0\\
       0&q^{-1}\end {array}\right ]
\left [\begin {array}{cc} \cos(\theta)&-\sin(\theta)\\
       \sin(\theta)&\cos(\theta)\end {array}\right ]
\left [\begin {array}{cc} {\it q'}^{-1}&0\\
      0&{\it q'}\end {array}\right ]
\ee
{}From (\ref{E:shift}), the action of $(t_l I)^2$ is a shift of
$\theta$ by $\pi$ and is indeed of order 2.  The parameters $q$ and
$q'$ are free for two of the matrices $[A,B,C]$, say $B$ and $C$.  If
$\delta(B)=\delta(C)$, then $\delta(A)=\delta'(C)$ and
$\delta'(A)=\delta'(B)$, and the composition law (\ref{E:slc})
coincides with the addition on $\theta$, and we have a special case of
solution (\ref{E:asym6}).  The solution (3.1) of \cite{Color}
corresponds to a slightly different splitting:
\be
\hat X =
\frac12\left [\begin {array}{cc} e^{-q}&e^q\\
       -e^{-q}&e^q\end {array}\right ]
\left [\begin {array}{cc} \cos(\theta)&-\sin(\theta)\\
       \sin(\theta)&\cos(\theta)\end {array}\right ]
\left [\begin {array}{cc} e^{q'}&-e^{q'}\\
      e^{-q'}&e^{-q'}\end {array}\right ]
\ee

As a summary we can state that any six vertex solution is one of the
two presented here, depending on whether $p_9$ vanishes or not.

\section{An eight-vertex Ansatz}
We take next a particular eight-vertex Ansatz: the three matrices
$A,B,C$ are assumed to be symmetric with respect to the secondary
diagonal:
\begin{equation}
\label{ansatz8}
	X=\begin{pmatrix}	X_{11}	&0	&0	&X_{14}	\cr
			0	&X_{22}	&X_{23}	&0	\cr
			0	&X_{32}	&X_{22}	&0	\cr
			X_{41}	&0	&0	&X_{11}
	                \cr	\end{pmatrix}
\end{equation}
There are many solutions  of (\ref{eq11}) having this form; here we
will analyze the solution that has the maximum number of parameters,
which is seven after the scalings have been fixed by putting $X_{11}=1$.

\subsection{General observations}
The solution $\S:=(A,B,C)$ can be given in terms of seven independent
parameters $a,\,b,\,c,$ $x,\,y,\,z,\,v$:
\begin{eqnarray} \displaystyle
\label{sol7a}
A&=&\left (
\begin {array}{cccc}
1&0&0&a\\0&x&{ \displaystyle\frac {b\left (v-x\right )}{ cy}}&0\\0&{
\displaystyle\frac {c\left (v-xyz\right )}{bz}}&x&0\\{
\displaystyle\frac {\left (v- y\right )\left (v-z\right )}{ayz}}&0&0&1
\end {array}\right ) \\
\label{sol7b}
B&=&\left (
\begin {array}{cccc} 1&0&0&b\\0&y&{\displaystyle\frac {a\left (v-x\right )}{
cx}}&0\\0&{\displaystyle\frac {c\left (v-z\right
)}{az}}&y&0\\{\displaystyle\frac {\left (v-y
\right )\left (v-xyz\right )}{bxz}}&0&0&1
\end {array}\right ) \\
\label{sol7c}
C&=&\left (
\begin {array}{cccc} 1&0&0&c\\0&z&{\displaystyle\frac {a\left (v-xyz\right
)}{bx}}&0\\0&{\displaystyle\frac {b\left (v-z\right
)}{ay}}&z&0\\{\displaystyle\frac {\left (v- y\right )\left (v-x\right
)}{cxy}}&0&0&1
\end {array}\right )
\end{eqnarray}

This solution is globally stable by the diagonal group $\g_d$, and the
action of $g\in\g_d $ (\ref{E:diagonal}) moves only $a,b,c$ as
follows:
\begin{eqnarray}
a \mapsto a\;  t_1^{-2} t_2^{-2}, \quad
b \mapsto b\;  t_1^{-2} t_3^{-2}, \quad
c \mapsto c\;  t_2^{-2} t_3^{-2}
\end{eqnarray}
The remaining four parameters are gauge invariant and are therefore
coordinates on the orbit space $\eta=\S/\g_d$.  The choice of a gauge
amounts to the choice of three functions $ a(x,y,z,v)$, $ b(x,y,z,v)$,
$ c(x,y,z,v)$, and we will later see the effect of choosing a
particular form.

The action of the generators of $\A ut$ is a birational transformation
of the parameters and reads:
\bea
K_a&:&[a,b,c] \mapsto [{-\displaystyle\frac {\left (v-y\right
 )\left (v-z\right )}{ayz}},{\displaystyle\frac {c\left ( v-z\right
 )}{az}},{\displaystyle\frac {b\left (v-z\right )}{ay}}]
\nonumber\\
&&[x,y,z,v]  \mapsto
[{\displaystyle\frac {\left (v-z-y
\right )x}{v-x-xyz}},y,z,z-v+y]
\label{ka}
\eea
\bea
K_b&:&[a,b,c] \mapsto [{\displaystyle\frac {c\left (v-xyz\right
 )}{bz}},-{\displaystyle\frac {\left (v-z-x\right )
\left (v-y\right )\left (v-xyz\right )}{bxz\left (v-xyz-y\right )}},{
\displaystyle\frac {a\left (v-xyz\right )}{bx}}]
\nonumber\\
\label{kb}
&&[x,y,z,v]  \mapsto
[x,{\displaystyle\frac {\left (v-z-x\right )y}{v-xy
z-y}},z,z-v+x]
\eea
\bea
K_c&:&[a,b,c] \mapsto [{\displaystyle\frac {b\left (v-x\right
 )}{cy}},{\displaystyle\frac {a\left (v-x\right )}{cx}}
 ,-{\displaystyle\frac {\left (v-y\right )\left (v-x\right )}{cxy}}]
\nonumber\\ \label{kc}
&&[x,y,z,v]  \mapsto
[x,y,{\displaystyle\frac {\left (v-x-y\right )z}{v-xyz-z}},x-v+y]
\eea
One verifies here that $\A ut$ indeed acts on the orbit space $\eta$,
since the transformation of $x$, $y$, $z$ and $v$  does {\it not} depend
on $a,b,c$.

There are two invariants of $\A ut$ on $\eta$:
\begin{equation}
 \D_1 = \frac{ v( 2v - xyz - x - y - z) }{xyz}
\label{E:del1}
\end{equation}
\begin{equation}
 \D_2 = \frac{(v-x)(v-y)(v-z)(v-xyz)}{x^2 y^2 z^2}
\label{E:del2}
\end{equation}
There is a canonical way to find these invariants~\cite{FaVi93}. It
consists of first calculating the squares of the generators $K_a, K_b,
K_c$ in homogeneous coordinates.  Such squares appear as the
multiplication by some polynomial ($\Phi_a, \Phi_b, \Phi_c$). For
example using the homogenizing variable $t$:
\begin{eqnarray*}
 K_a:	& t &	\mapsto	t\left (t^{2}x-t^{2}v+xyz\right ) \\
	& x &	\mapsto	xt^{2}\left (z-v+y\right ) \\
	& y &	\mapsto	y\left (t^{2}x-t^{2}v+xyz\right ) \\
	& z &	\mapsto	z\left (t^{2}x-t^{2}v+xyz\right ) \\
  & v &	\mapsto	\left (z-v+y\right )\left (t^{2}x-t^{2}v+xyz\right ) \\
K_a^2 \simeq \Phi_a & = &
t^{4}\left (t^{2}x-t^{2}v+xyz\right )^{3}v\left (z-v+y\right )
\end{eqnarray*}
Any rational invariant is the ratio of two polynomials which have the
same covariance properties under $K_a$ (resp.\ $K_b$ and $K_c$).
The covariance factors are known to be the factors appearing in
$\Phi_a$, (resp.\ $\Phi_b$, $\Phi_c$), and for a given degree, there
are only a finite number of possible covariance factors. It is thus
possible to find all algebraic invariants of a given degree. This
algorithm is unfortunately unbounded, since we do  not know any
bound on the  degree of the invariant. However, it proves quite efficient
in practice.

The two invariants are `global', as they can be calculated from any of
the three matrices $A,B$ or $C$ using the polynomials $p_5$ and $p_9$
as given in \cite{BeMaVi92}:
\[
\D_1=  - 2\, \frac{ p_9}{p_5}, \quad
\D_2= \frac{\mbox{product of anti-diagonal entries}}
  {\mbox{product of diagonal entries}}.
\]
The surfaces $\D_1=\mbox{constant}$, $\D_2=\mbox{constant}$, in $\eta$
are preserved by the induced action of $\A ut$.  They are of generic
dimension 2 and define the varieties where the spectral parameters
live.  The invariants $ \D_1$
and $ \D_2$ are the coordinates on the moduli space of the solution
$\S$. Note that they are invariant by any permutation of $x,y,z$.
Note also that the free fermion condition~\cite{FaWu70} is just
$\D_1=0$.

For later discussions let us introduce the parameters $q_i$ by
\be
q_1=\frac{v-x}{\sqrt{xyz}},\,
q_2=\frac{v-xyz}{\sqrt{xyz}},\,
q_3=\frac{v-z}{\sqrt{xyz}},\,
q_4=\frac{v-y}{\sqrt{xyz}}.
\label{E:qdef}
\ee
For the inverse relations define first $\Lambda$ by
\be
\prod_{i=1}^4 (\Lambda-q_i)=1,
\ee
and then
\be
x=(\Lambda-q_1)(\Lambda-q_2),\,
y=(\Lambda-q_4)(\Lambda-q_2),\,
z=(\Lambda-q_3)(\Lambda-q_2),\,
v=\Lambda(\Lambda-q_2),
\ee
and for the $\Delta$'s we get
\bea
\Delta_2&=&q_1q_2q_3q_4,\\
\Delta_1&=&-2\Lambda^2+\Lambda(q_1+q_2+q_3+q_4)\nonumber\\
\label{E:d1q}&=&
q_1q_2+q_3q_4-x-\frac1x=q_1q_3+q_2q_4-y-\frac1y=
q_2q_3+q_1q_4-z-\frac1z.
\eea
With these definitions the anti-diagonal entries of our solution can
be written as
\bea
ad(A)&=&\left\{a,\frac{b}{c}\sqrt{\frac{xz}{y}}\, q_1,
\frac{c}{b}\sqrt{\frac{xy}{z}}\, q_2,\frac{x}{a}\, q_3 q_4\right\},\\
ad(B)&=&\left\{b,\frac{a}{c}\sqrt{\frac{yz}{x}}\, q_1,
\frac{c}{a}\sqrt{\frac{xy}{z}}\, q_3,\frac{y}{b}\, q_2 q_4\right\},\\
ad(C)&=&\left\{c,\frac{a}{b}\sqrt{\frac{yz}{x}}\, q_2,
\frac{b}{a}\sqrt{\frac{xz}{y}}\, q_3,\frac{z}{c}\, q_1 q_4\right\}.
\eea
Note that the $q_i$ behave almost like the color parameters.

\subsection{Specific gauges and related parametrizations}
We will show here how the gauge condition, i.e. the choice of $a,b,c$
as functions of $x,y,z,v$ affects the distribution of the parameters
among the three members of the solution triplet. The solution we have
is a four parameter solution, once the gauge is fixed, as is Baxter's
solution~\cite{Ba81,Ba71}.  Note that for us, `fixing the gauge' means
preventing continuous residual gauge freedom but leaves room for
discrete transformations.

As a possible simple gauge we could take $a=1, b=1, c=1$.  With this
choice, choosing $C$ uses up {\it all} four parameters of the
solution, and $B$ and $A$ are completely determined once $C$ is known,
thus the parameter content in this case is (4,0,0). The solution is
fully rational but does not lead to a parametrized family of commuting
transfer matrices. It leads to another very interesting --- although
apparently less constrained --- situation: we have an infinite
sequence of commuting transfer matrices with commutation between
successors. We shall not explore this possibility here.

An important question now is the following: How should we choose the
gauge condition in order to get a parametrized family of solutions?
Clearly the minimum requirement is to choose $a,b,c$ in such a way
that the knowledge of $C$, for example, uses only {\em three} of the
four available parameters on $\eta$, leading to the parameter content
(3,1,0).  In other words, the gauge choice must lower the rank of the
set of anti-diagonal elements
\be
\Sigma:=\left\{z,c,\frac{a}{b}\sqrt{\frac{yz}{x}}\, q_2,
\frac{b}{a}\sqrt{\frac{xz}{y}}\, q_3,\frac{z}{c}\, q_1 q_4\right\},
\ee
from four to three.  Instead of $\Sigma$ we could consider the set
$\Sigma':=\{z,\, c,\,\frac{b}{a}\sqrt{\frac{x}{y}}\, q_3,\,
q_2q_3,\,q_1q_4\}$, or using (\ref{E:d1q}), $\Sigma'' := \{z,\, c,\,
\frac{b}{a}\sqrt{\frac{x}{y}}\, q_3,\, \Delta_1,\, \Delta_2\}$.
In this last set, $z$, $\D_1$, and $\D_2$ are clearly functionally
independent, so in order to have no more than these three parameters we
must impose the condition
\be
\label{eqc}
c= f(z,\D_1,\D_2) ,\quad
\frac{a}{b} = \sqrt{\frac{x}{y}}\, q_3 \,  g(z,\D_1,\D_2),
\ee
where $f$ and $g$ are some arbitrary functions.

It was argued earlier that for many applications $C$ and $B$, say,
should have a similar structure, and in particular the same number of
free parameters. If we therefore apply the above argument to $B$ we
get in a similar way the conditions
\be
\label{eqb}
b = h(y, \D_1,\D_2) , \quad
\frac{a}{c} = \sqrt{\frac{x}{z}}\, q_3 \,  k(y, \D_1,\D_2),
\ee
where $h$ and $k$ are free functions.  The compatibility of
(\ref{eqc}, \ref{eqb}) (solving for $a$ in two ways) implies
\begin{equation}
k =\omega(\D_1,\D_2)  h(y,\D_1,\D_2)/\sqrt{y}, \quad
g = \omega(\D_1,\D_2) f(z,\D_1,\D_2)/ \sqrt{z},
\end{equation}
where $\omega$ is an arbitrary function, so that
\begin{equation}
a =\sqrt{\frac{x}{yz}}\;q_3\; h(y)\;f(z)\;\omega
=\frac{v-z}{yz}\; h(y)\;f(z)\;\omega.
\end{equation}
[From now on we do not write out the $\Delta_1,\Delta_2$ dependence.]

In order to get a true one parameter family of commuting transfer
matrices we want the matrices $B$ and $C$ to be in the same
parametrized family, i.e: $B=R(y,\D_1, \D_2)$ and $C=R(z,\D_1, \D_2)$
for some $R(\tau, \D_1, \D_2)$. This can be done, the condition is
$h=f$ and yields:
\be
R(\tau, \D_1, \D_2)
 = \begin{pmatrix} 1& 0& 0& f(\tau) \cr
	0& \tau\vphantom{\dfrac1T}& \omega\, f(\tau)\,Q & 0 \cr
	0& \dfrac{\tau }{\omega\,f(\tau)} & \tau & 0 \cr
	\dfrac{\tau\Delta_2}{f(\tau)\,Q}& 0& 0& 1 \end{pmatrix}
\ee
with $Q$ a root of
\be
\tau Q^2-\tau^2 Q-\tau  Q\Delta_1-Q+\tau\Delta_2 =0.
\label{E:ellc8v}
\ee

This result shows that the elliptic curve (\ref{E:ellc8v}) must be
introduced even if we just want to write two entries of the solution
as a parametrized family.  This still leaves considerable freedom in
choosing the gauge. Simple-looking results are obtained e.g.\ if we
take $\omega=1$ and $f=1$ or $\tau$, but these are no longer rational
in $\tau$ since $Q(\tau)$ will involve square roots.

The above expressions for $a,\,b,\,c$, were obtained by the condition
that $B$ and $C$ depend only on three of the four parameters. If we now
continue and require the same on $A$ we obtain the additional
conditions
\begin{equation}
a=m(x,\Delta_1,\Delta_2),\quad \frac{c}b=\sqrt{\frac{z}{y}}\;q_1\;
n(x,\Delta_1,\Delta_2).
\end{equation}
The compatibility of these with (\ref{eqc},\ref{eqb}) leads eventually
to the symmetric solution
\begin{equation}
a=\sqrt{\frac{x}{q_1q_2}}\;\phi_1\phi_2,\quad
b=\sqrt{\frac{y}{q_1q_3}}\;\phi_1\phi_3,\quad
c=\sqrt{\frac{z}{q_2q_3}}\;\phi_2\phi_3,\quad
\label{E:fg}
\end{equation}
where $\phi_i=\phi_i(\Delta_1,\Delta_2)$ is the residual freedom in
the choice of gauge and the previously used function $\omega$ is related
to $\phi_3$ by $\omega\phi_3^2=1$.

The gauge choice (\ref{E:fg}) leads to the following parameter
counting: each matrix of the triplet contains three independent
entries, and the choice of $C$, say, fixes three of the four
parameters of the solution. One free parameter is left for $B$ and
finally $A$ is determined once $B$ is chosen.

\subsection{Elliptic parametrization}
Baxter solution~\cite{Ba71} is actually contained in (\ref{sol7a},
\ref{sol7b}, \ref{sol7c}) if we choose the gauge (\ref{E:fg}) with
$\phi_i=\Delta_2^{1/4}$, in other words
$f(\tau)=h(\tau)=m(\tau)=\sqrt{\tau\Delta_2/Q(\tau)}$ and
$\omega=1/\sqrt{\Delta_2}$. This gauge is also uniquely defined by the
requirement that the matrices are symmetric under the usual
transposition.

But fixing the gauge is not the end of the story. For a good spectral
parameter we need also a good composition rule, in this case it is
obtained as follows. Fixing $C$ fixes the values of the invariants
$\D_1(C)$ and $\D_2(C)$, and therefore the elliptic curve
(\ref{E:ellc8v}). $C$, $B$, and $A$ will then be given by three points on
this same curve. There is a natural addition rule on elliptic curves,
and to verify that our parameters satisfy it we have to use the usual
uniformization with elliptic functions Baxter~\cite{Ba71}. The result
is as follows: Let us define $\gamma$ and $k$ by
\begin{equation}
\label{ellipd1d2}
\D_1 = - 2 cn(\gamma) dn(\gamma),  \quad \D_2= sn^4(\gamma) k^2
\end{equation}
where $sn, cn, dn$ are the Jacobi elliptic functions of modulus $k$,
and $\sigma,\,\rho,\chi$ by
\begin{equation}
\label{ellipyz}
z= \frac{sn(\sigma)}{ sn(\gamma - \sigma)},\quad
y= \frac{ sn(\rho)}{ sn(\gamma - \rho)}, \quad
x= \frac{ sn(\chi)}{ sn(\gamma - \chi)}, \quad
\end{equation}
then the relations (\ref{E:del1},\ref{E:del2}) are satisfied when we
take
\be
\label{ellipv}
v=\frac{sn(\rho)[sn(\sigma)sn(\chi)+sn(\gamma)sn(\gamma-\rho)]}{
sn(\gamma-\chi)sn(\gamma-\rho)sn(\gamma-\sigma)}
\ee
and use the addition rule
\be
\rho=\sigma+\chi.
\label{E:spe}
\ee
To complete the parametrization we note that
\bea
q_1&=&sn(\gamma)\sqrt{\frac{sn(\sigma)sn(\gamma-\sigma)}{
sn(\rho)sn(\gamma-\rho)sn(\chi)sn(\gamma-\chi)}}\,,\\
q_2&=&sn(\gamma)\sqrt{\frac{sn(\rho)sn(\gamma-\rho)}{
sn(\sigma)sn(\gamma-\sigma)sn(\chi)sn(\gamma-\chi)}}\,,\\
q_3&=&sn(\gamma)\sqrt{\frac{sn(\chi)sn(\gamma-\chi)}{
sn(\sigma)sn(\gamma-\sigma)sn(\rho)sn(\gamma-\rho)}}\,,\\
q_4&=&sn(\gamma)k^2\sqrt{sn(\chi)sn(\gamma-\chi)
sn(\rho)sn(\gamma-\rho)sn(\sigma)sn(\gamma-\sigma)}\,,
\eea
\be
Q(x)=\frac{sn(\gamma)^2}{sn(\chi)sn(\gamma-\chi)}\,,
\ee
and similarly for $y,\,z$.

If we now define (note the overall scaling)
\be
R(\alpha,\gamma,k)=
 \begin{pmatrix} sn(\gamma-\alpha)& 0& 0& sn(\alpha)\,sn(\gamma)\,k \cr
	0& sn(\alpha) & sn(\gamma) & 0 \cr
	0& sn(\gamma) & sn(\alpha) & 0 \cr
	sn(\alpha)\,sn(\gamma)\,k & 0& 0& sn(\gamma-\alpha)& \end{pmatrix}
\ee
then $A= R(\chi,\gamma,k), B= R(\rho,\gamma,k), C= R(\sigma,\gamma,k)$
solve (\ref{E:abcsh}), and this is exactly Baxter's solution.  [If $k=0$
we get a special case of the six-vertex solution.]

What (\ref{ellipd1d2}), (\ref{ellipyz}), (\ref{ellipv}), together with
(\ref{E:spe}) show is that the two dimensional surface in $\eta$ given
by fixing the values of $\Delta_1$ and $\Delta_2$ is a product of two
elliptic curves (or two points on the same elliptic curve). Baxter's
parametrization makes it explicit.  Furthermore it also allows to
visualize the action of $\A ut$, since
\be
\begin{array}{ll}
K_a K_b: \rho \mapsto \rho + \gamma, &\sigma \mapsto \sigma \\
K_b K_c:  \rho \mapsto \rho ,&\sigma \mapsto \sigma + \gamma
\end{array}
\ee

\section{Conclusion}
We have analyzed several two-state solutions of the Yang-Baxter
equations, and shown how starting from a rational solution without
recognizable structure one can construct spectral and moduli
parameters using the symmetries of the equations. For five- and
six-vertex ansatze our results are complete.

The effect of gauge transformation on the parametrization is
particularly interesting.  One of the lessons of the present work is
that it can be much easier to solve the YBE when a gauge has not been
fixed. Finding a good parametrization is a separate problem, which can
be done at leisure, after a solution has been found.

The present method can be used for any solution of the Yang-Baxter
equations, whenever they are found. Unfortunately we do not yet have a
thorough analysis of these equations with absolutely no a priori
assumptions on their form (i.e. no Ansatz at all), in the spirit of
the complete resolution of the `constant' equations obtained earlier
by one of the authors~\cite{Hi92}.

\subsection*{Acknowledgments} This work was started when one of the
authors (JH) was visiting L.P.T.H.E.  under the support of CNRS and
Academy of Finland. On of us (CV) would like to thank J-M. Maillard
for discussions and P.A. Pearce for pointing out ref(\cite{No92}).

\label{lastpage}
\end{document}